\documentclass[preprint,12pt]{elsarticle}

\usepackage{amssymb}
\usepackage{amsthm}
\usepackage{amsmath}
\usepackage{fancybox}

\newtheorem{thm}{Theorem}
\newtheorem{lem}{Lemma}
\newdefinition{rem}{Remark}
\newdefinition{exa}{Example}
\newdefinition{defn}{Definition}
\newproof{pf}{{\bf Proof}}
\newproof{pthone}{{\bf Proof of Theorem~\protect\ref{th:1}}}
\newproof{pthscnd}{{\bf Proof of Theorem~\protect\ref{th:2}}}
\newcommand{\eq}[1]{\mbox{\rm(\ref{#1})}}
\newcommand{\rra}{\rightrightarrows}

\catcode`\@=11  \@addtoreset{equation}{section}  \catcode`\@=12
\newcommand{\RB}{\mathbb{R}}   \newcommand{\NB}{\mathbb{N}}
\newcommand{\QB}{\mathbb{Q}}
\newcommand{\DC}{\mathcal{D}}  \newcommand{\PC}{\mathcal{P}}
\newcommand{\UC}{\mathcal{U}}
\newcommand{\vep}{\varepsilon}
\newcommand{\vfi}{\varphi}

\journal{arXiv\quad}

\begin{document}

\begin{frontmatter}

\title{The approximate variation of univariate uniform space\\
        valued functions and pointwise selection principles}

\author[hsenn]{Vyacheslav V.~Chistyakov\corref{cor1}}
\ead{czeslaw@mail.ru, vchistyakov@hse.ru}

\author[hsenn]{Svetlana A.~Chistyakova}
\ead{schistyakova@hse.ru}


\cortext[cor1]{Corresponding author.}

\address[hsenn]{Department of Informatics, Mathematics, and Computer Science,\\
National Research University Higher School of Economics,\\
Bol'shaya Pech{\"e}rskaya Street 25/12, Nizhny Novgorod\\
603155, Russian Federation}

\begin{abstract}
Let $T\subset\mathbb{R}$ and $(X,\mathcal{U})$ be a uniform
space with an at most countable gage of pseudometrics
$\{d_p:p\in\mathcal{P}\}$ of the uniformity $\mathcal{U}$.
Given $f\in X^T$ (=\,the family of all functions from $T$ into $X$),
the {\em approximate variation\/} of $f$ is the two-parameter
family $\{V_{\varepsilon,p}(f):\varepsilon>0,p\in\mathcal{P}\}$,
where $V_{\varepsilon,p}(f)$ is the greatest lower bound of
Jordan's variations $V_p(g)$ on $T$ with respect to $d_p$ of all
functions $g\in X^T$ such that $d_p(f(t),g(t))\le\varepsilon$
for all $t\in T$. We establish the following pointwise selection principle:
{\em If a pointwise relatively sequentially compact sequence
of functions $\{f_j\}_{j=1}^\infty\subset X^T$ is such that
$\limsup_{j\to\infty}V_{\varepsilon,p}(f_j)<\infty$ for all
$\varepsilon>0$ and $p\in\PC$, then it contains a subsequence
which converges pointwise on $T$ to a bounded regulated
function $f\in X^T$.} We illustrate this result by appropriate
examples, and present a characterization of regulated
functions $f\in X^T$ in terms of the approximate variation.
\end{abstract}

\begin{keyword}
uniform space \sep approximate variation \sep  regulated function \sep
pointwise convergence \sep selection principle

{\em MSC\,2010:} 26A45 \sep 54E15 \sep 40A30 \sep 26A30
\end{keyword}
\end{frontmatter}

\section{Introduction}

The historically first extensions of Bolzano-Weierstrass' theorem (viz., a bounded
sequence in $\RB$ admits a convergent subsequence) are two theorems of Helly
\cite{Helly} (cf.\ also \cite[II.8.9--10]{Hild}, \cite[VIII.4.2]{Nat}) for real monotone
functions and for functions of bounded (Jordan) variation on $I=[a,b]$,
conventionally called Helly's selection principles.%
\footnote{We recall Helly's selection principle for monotone functions
in Section~\ref{s:proofs}.}
Their significant role in analysis and topology is well-known
(\cite{Dieu}, \cite{Hild},  \cite{Kelley}, \cite{Nat}, \cite{Schwartz}). 
A vast literature already exists concerning generalizations of Helly's
principles for various classes of functions (\cite{MatSb}, \cite{Sovae}, \cite{JMAA05},
\cite{MMS}, \cite{Studia}, \cite{JMAA17}, \cite{Wat}, \cite{Fr}, \cite{MuOr}, and
references therein, mostly for metric space valued functions) as well as their applications
in the theory of convergence of Fourier series, stochastic processes, Riemann- and
Lebesgue-Stieltjes integrals, optimization, set-valued analysis, generalized ODEs,
modular analysis (\cite{Aumann}, \cite{Barbu}, \cite{Chant}, \cite{Sovae}, 
\cite{MMS}, \cite{JMAA19}, \cite{JMAA07}, \cite{Jeff}, \cite{Schwabik}). 

\smallbreak
The latter series of references above exhibits also an important role played by
regulated functions in various areas of mathematics.
Classically, a function $f:I=[a,b]\to\RB$ is said to be
regulated provided the left limit $f(t-0)\in\RB$ exists at each point $a<t\le b$ and
the right limit $f(t+0)\in\RB$ exists at each point $a\le t<b$. It is well known that
each regulated function on $I$ is bounded, has at most a countable set of
discontinuity points, and is the uniform limit of a sequence of step functions on~$I$.
Furthermore, there are different descriptions of regulated real and metric space
valued functions (\cite{Aumann}, \cite{Chant}, \cite{JMAA05}, \cite{DAN06}, \cite{MMS},
\cite{Studia}, \cite{Wat}, \cite{Dieu}, \cite{Dudley}, \cite{Fr}, \cite{Goff},
\cite{Jeff}, \cite{MuOr}, \cite{Perl}, \cite{Tolsto}, \cite{Tretya}).%
\footnote{In Section~\ref{s:mains}, we adopt a more general definition of
regulated functions.}

\smallbreak
There is an explicit indication in the literature (e.g., \cite{arXiv19}, \cite{Wat}) of a certain
conjunction of `characterizations of regulated functions' and `pointwise selection
principles'. It can be formulated (heuristically, at present) as follows: if we have
a characterization of regulated functions in certain {\em terms}, then, in the same
{\em terms}, we may obtain a (Helly-type) pointwise selection principle, possibly
outside the scope of regulated functions. This is demonstrated by the main results
of this paper, Theorems~\ref{th:1} and~\ref{th:2}, as well.

\smallbreak
The class of uniform spaces is `situated' properly between the class of metric spaces
and the class of topological spaces (\cite{Bour}, \cite{Kelley}). Applications of uniform
spaces outside topology seem to be not numerous (e.g., \cite{DAN06}, \cite{MatTr},
\cite{Tolsto}, \cite{Tretya}), and the most important of them
being to topological groups (\cite{Arhan}, \cite{Bour}). However, since uniform
spaces can be described in terms of their gages of pseudometrics (see \cite{Kelley} and
Section~\ref{s:mains}), they form natural codomains (=\,ranges) for regulated
functions  and sequences of functions in selection principles. In fact, we may
define oscillations of functions, uniform pseudometrics,
bounded functions, functions of bounded (Jordan) variation and regulated functions
with respect to pseudometrics from the gage, pointwise and uniform convergences
of sequences of functions valued in a uniform space.

\smallbreak
The purpose of this paper is to present a new characterization of a uniform space
valued univariate regulated functions and, on the basis of it, establish a pointwise selection
principle in terms of the (new notion of) approximate variation for uniform space valued
sequences of univariate functions.
Further comments on the novelty of our results
can be found in Section~\ref{s:mains}.

\smallbreak
The paper is organized as follows. In Section~\ref{s:mains}, we present necessary
definitions and our main results, Theorems~\ref{th:1} and~\ref{th:2}.
In Section~\ref{s:prappvar}, we establish essential properties of the approximate variation
(Lemmas \ref{l:ess} and \ref{l:unif}), prove Theorem~\ref{th:1}, and present illustrating
examples including Dirichlet's function. In Section~\ref{s:proofs}, we prove
Theorem~\ref{th:2}, show that its main assumption \eq{e:Vepfj} is necessary for the
uniform convergence, but not necessary for the pointwise convergence of sequences of
functions, and comment on the applicability of Theorem~\ref{th:2} to sequences of
non-regulated functions.

\section{Definitions and main results} \label{s:mains}

We begin by reviewing certain definitions and facts needed for our results.

\smallbreak
Throughout the paper we assume (if not stated otherwise) that the pair
$(X,\UC)$ is a Hausdorff uniform space with the gage of pseudometrics
$\{d_p\}_{p\in\PC}$ (on $X$) of the uniformity $\UC$ for some index set $\PC$
(cf. \cite[Chapter~6]{Kelley}). For the reader's convenience, we describe this
assumption in more detail (I)--(IV).

\smallbreak
(I) Given $p\in\PC$, the function $d_p:X\times X\to[0,\infty)$ is a {\em pseudometric\/}
on $X$, i.e., for all $x,y,z\in X$, it satisfies the three conditions: $d_p(x,x)=0$,
$d_p(y,x)=d_p(x,y)$, and $d_p(x,y)\le d_p(x,z)+d_p(z,y)$ (triangle inequality).

\smallbreak
(II) The family $\{d_p\}_{p\in\PC}$ of pseudometrics on $X$ is the {\em gage\/} of the
uniformity $\UC$ for $X$ if $\{d_p\}_{p\in\PC}$ is the family of {\em all\/}
pseudometrics which are uniformly continuous on $X\times X$ relative to the product
uniformity derived from~$\UC$. Recall (\cite[Theorem~6.11]{Kelley}) that a pseudometric
$d$ on $X$ is {\em uniformly continuous on $X\times X$ relative to the product
uniformity\/} if and only if the set $U(d,\vep)=\{(x,y)\in X\times X:d(x,y)<\vep\}$
belongs to $\UC$ for all $\vep>0$. The uniform space $(X,\UC)$ is said to be
{\em Hausdorff\/} if conditions $x,y\in X$ and $d_p(x,y)=0$ for all $p\in\PC$
imply $x=y$.

\smallbreak
(III) Every family $\mathcal{F}$ of pseudometrics on $X$ generates a uniformity $\UC$
(by \cite[Theorem 6.15]{Kelley}). A direct description of the gage $\{d_p\}_{p\in\PC}$
of $\UC$ generated by $\mathcal{F}$ is as follows (\cite[Theorem 6.18]{Kelley}):
the family $\{U(d,\vep):\mbox{$d\in\mathcal{F}$ and $\vep>0$}\}$ is a subbase for
the uniformity~$\UC$; thus, given a pseudometric $d$ on $X$, $d\in\{d_p\}_{p\in\PC}$
if and only if, for each $\vep>0$, there are $\delta>0$, $n\in\NB$ and
$p_1,\dots,p_n\in\PC$ such that $\bigcap_{i=1}^n U(d_{p_i},\delta)\subset U(d,\vep)$.
Furthermore (\cite[Theorem 6.19]{Kelley}), if $\{d_p\}_{p\in\PC}$ is the gage of $\UC$,
then the family $\{U(d_p,\vep):\mbox{$p\in\PC$ and $\vep>0$}\}$ is a base for the
uniformity $\UC$ (i.e., for every $U\in\UC$ there are $p\in\PC$ and $\vep>0$ such that
$U(d_p,\vep)\subset U$).

\smallbreak
(IV) ``Each concept which is based on the notion of a uniformity can be described
in terms of a gage because each uniformity is completely determined by its gage''
(cf. \cite[Theorem~6.19]{Kelley}). 

\smallbreak
A particular case of a uniform space is a metric space $(X,d)$ with metric~$d$.
The standard base of the uniformity $\UC$ for $X$ is the family
$\{U(d,\vep):\vep>0\}$, where $U(d,\vep)$ is as in (II). The uniformity $\UC$ has a
countable base, e.g., the family $\{U(d,{1/k}):k\in\NB\}$. Two more examples are
in order. First, the family $\{U_\alpha:\alpha\in\RB\}$ with
$U_\alpha=\{(x,y)\in\RB\times\RB:\mbox{$x=y$ or $x>\alpha$, $y>\alpha$}\}$
is a base of a uniformity on $\RB$. Second, let $\mbox{\rm CB}(X;\RB)$ be the family
of all continuous bounded functions from a completely regular topological space $X$
into $\RB$ and, given $p\in\PC\equiv\mbox{\rm CB}(X;\RB)$, define
$d_p(x,y)=|p(x)-p(y)|$ for $x,y\in X$. The family $\{d_p:p\in\PC\}$ generates
a gage for a uniformity on $X$ which is compatible with the topology on~$X$.
For more examples, we refer to \cite{Bour,Engel,Kelley,Will}.

\smallbreak
For the sake of brevity, we often write $X$ in place of $(X,\UC)$.

\smallbreak
Given a nonempty set $T$ (in what follows, $T\subset\RB$), we denote by $X^T$ the
family of all functions $f:T\to X$ mapping $T$ into $X$ and equip it with (extended-valued)
{\em uniform pseudometrics}
  $$d_{T,\,p}(f,g)=\sup_{t\in T}d_p(f(t),g(t)),\quad\,\, f,g\in X^T,\quad p\in\PC.$$

\smallbreak
For $p\in\PC$, the {\em oscillation of $f\in X^T$ with respect to pseudometric $d_p$}
is the quantity
  $$|f(T)|_p=\sup_{s,t\in T}d_p(f(s),f(t))\in[0,\infty],$$
also known as the {\em$d_p$-diameter\/} of the image $f(T)=\{f(t):t\in T\}\subset X$.
We denote by $\mbox{\rm B}_p(T;X)=\{f\in X^T:|f(T)|_p<\infty\}$ the family of all
{\em$d_p$-boun\-ded functions\/}, and by
$\mbox{\rm B}(T;X)=\bigcap_{p\in\PC}\mbox{\rm B}_p(T;X)$---the family of all
{\em bounded functions\/} mapping $T$ into~$X$.

\smallbreak
Recall that a sequence of points $\{x_j\}\equiv\{x_j\}_{j=1}^\infty\subset X$
{\em converges\/} in the uniform space $X$ to an element $x\in X$ (as $j\to\infty$)
if and only if $d_p(x_j,x)\to0$ as $j\to\infty$ for all $p\in\PC$. Since $X$ is Hausdorff,
the limit $x$ is determined uniquely.
A subset $Y$ of the uniform space $X$ is said to be {\em sequentially compact\/}
({\em relatively sequentially compact}) if each sequence of points from $Y$ has a
subsequence which converges in $X$ to an element of $Y$ (to an element of $X$,
respectively).

\smallbreak
Suppose we have a sequence of functions $\{f_j\}\equiv\{f_j\}_{j=1}^\infty\subset X^T$
and $f\in X^T$. If $p\in\PC$, we write: (a) {\em $f_j\to_p f$ on~$T$} if
$d_p(f_j(t),f(t))\to0$ as $j\to\infty$ for all $t\in T$, and (b) {\em $f_j\rra_p f$ on $T$} if
$d_{T,\,p}(f_j,f)\to0$ as $j\to\infty$. We denote by: (c) {\em $f_j\to f$ on~$T$} the
{\em pointwise\/} (or {\em everywhere\/}) {\em convergence\/} of $\{f_j\}$ to~$f$,
i.e., $f_j\to_p f$ on $T$ for all $p\in\PC$, and (d) {\em$f_j\rra f$ on $T$} the
{\em uniform convergence\/} of $\{f_j\}$ to~$f$, i.e., $f_j\rra_p f$ on $T$ for all
$p\in\PC$. Clearly, (d) implies (c), but not vice versa.

\smallbreak
A sequence of functions $\{f_j\}\subset X^T$ is said to be {\em pointwise relatively
sequentially compact\/} on $T$ if the set $Y=\{f_j(t):j\in\NB\}$ is relatively sequentially
compact in $X$ for all $t\in T$.

\smallbreak
From now on, we assume that $T$ is a nonempty subset of the reals~$\RB$.

\smallbreak
Given $f\in X^T$ and $p\in\PC$, the ({\em Jordan-type\/}) {\em variation of $f$
with respect to pseudometric $d_p$\/} is the quantity (e.g., \cite{MatSb}, \cite{MatTr},
\cite{Schwartz}, \cite{Tolsto})
  $$V_p(f,T)=\sup_{\pi}\sum_{i=1}^m d_p(f(t_i),f(t_{i-1}))\in[0,\infty],$$
where the supremum is taken over all partitions $\pi$ of $T$, i.e., $m\in\NB$ and
$\pi=\{t_i\}_{i=0}^m\subset T$ such that $t_{i-1}\le t_i$ for all $i=1,2,\dots,m$.
We denote by $\mbox{\rm BV}_p(T;X)=\{f\in X^T:V_p(f,T)<\infty\}$ the family of all
{\em functions\/} from $T$ into $X$ {\em of bounded $d_p$-variation}, and by
$\mbox{\rm BV}(T;X)=\bigcap_{p\in\PC}\mbox{\rm BV}_p(T;X)$---the family of all
{\em functions of bounded variation\/} (or BV functions, for short).

\smallbreak
The following definition is crucial for the whole subsequent material.

\begin{defn}
The {\em approximate variation\/} of a function $f\in X^T$ is the two-parameter
family $\{V_{\vep,\,p}(f,T):\mbox{$\vep>0$ and $p\in\PC$}\}\subset[0,\infty]$ of
{\em $\vep$-$d_p$-variations\/} $V_{\vep,\,p}(f,T)$ of $f$, given for $\vep>0$
and $p\in\PC$, by
  \begin{equation} \label{d:Vep}
V_{\vep,\,p}(f,T)=\inf\,\bigl\{V_p(g,T):\mbox{$g\in\mbox{\rm BV}_p(T;X)$
and $d_{T,\,p}(f,g)\le\vep$}\bigr\}
  \end{equation}
with the convention that $\inf\varnothing=\infty$.
\end{defn}

The notion of {\em$\vep$-variation\/} $V_\vep(f,T)$ of a function $f\in X^T$ with
$T=[a,b]$ a closed interval in $\RB$ and $X=\RB^N$ was introduced by
Fra{\v n}kov{\'a} \cite[Definition~3.2]{Fr}. For any $T\subset\RB$ and a metric space
$(X,d)$, this notion was considered and extended by the authors in
\cite[Sections 4 and 6]{Studia}, and the corresponding approximate variation was
thoroughly studied by the first author in \cite{arXiv19} (in particular, Example 3.10
from \cite{arXiv19} shows that the infimum in \eq{d:Vep} is in general not attained).
The notion of $\vep$-variation was also generalized by the authors for metric space
valued bivariate functions in \cite{JMAA17}.

\smallbreak
The notion of the approximate variation is useful in characterizing regulated metric space
valued functions in the usual sense (\cite{arXiv19}, \cite{Studia}, \cite{Fr}).
For our purposes, we adopt the following more general definition.

\smallbreak
Given $p\in\PC$, a function $f\in X^T$ is said to be {\em$d_p$-regulated on~$T$} (in
symbols, $f\in\mbox{\rm Reg}_p(T;X)$) if it satisfies the $d_p$-Cauchy conditions at
every left limit point of $T$ and every right limit point of~$T$. More explicitly,
given $\tau\in T$, which is a {\em left limit point for $T$} (i.e.,
$T\cap(\tau-\delta,\tau)\ne\varnothing$ for all $\delta>0$), we have
$d_p(f(s),f(t))\to0$ as $T\ni s,t\to\tau-0$, and similarly, given $\tau'\in T$,
which is a {\em right limit point for $T$} (i.e.,
$T\cap(\tau',\tau'+\delta)\ne\varnothing$ for all $\delta>0$), we have
$d_p(f(s),f(t))\to0$ as $T\ni s,t\to\tau'+0$. A function $f\in X^T$ is said to be
{\em regulated on~$T$} if it is $d_p$-regulated on $T$ for all $p\in\PC$.
We set $\mbox{\rm Reg}(T;X)=\bigcap_{p\in\PC}\mbox{\rm Reg}_p(T;X)$.
It is to be noted that a regulated function need not be bounded
in general (for instance, the function $f\in\RB^T$, given on the set
$T=[0,1]\cup\{2-(1/n):n=2,3,\dots\}$
by $f(t)=t$ for $0\le t\le1$ and $f(2-(1/n))=n$ for $n=2,3,\dots$, is regulated
on $T$ in the above sense, but not bounded). In the case when $T=I=[a,b]$ is a closed
interval in~$\RB$ and $X$ is complete we have: $f\in\mbox{\rm Reg}(I;X)$ if and only if
the (left) limit of $f(t)$ as $t\to\tau-0$ exists in $X$ at each point $\tau\in(a,b]$ 
and the (right) limit of $f(t)$ as $t\to\tau'+0$ exists in $X$ at each point
$\tau'\in[a,b)$ (cf.~\cite[Section~4]{MatTr}).

\smallbreak
Our first result is a characterization of regulated functions in terms of the approximate
variation (containing similar results from \cite[equality~(4.2)]{Studia}
and \cite[Proposition~3.4]{Fr} as particular cases):

\begin{thm} \label{th:1}
Suppose $\varnothing\ne T\subset\RB$ and $(X,\UC)$ is a Hausdorff uniform space.
Given $p\in\PC$, $\{f\!\in\! X^T\!:V_{\vep,\,p}(f,T)\!<\!\infty\,\,\forall\,%
\vep\!>\!0\}\subset\mbox{\rm Reg}_p(T;X)$, and this inclusion turns into the equality
for $T=I=[a,b]$. Consequently,
  \begin{equation*}
\mbox{\rm Reg}(I;X)=\{f\!\in\! X^I\!:\mbox{$V_{\vep,\,p}(f,I)\!<\!\infty$ for all
$\vep\!>\!0$ and $p\in\PC$}\}.
  \end{equation*}
\end{thm}

\smallbreak
As it was already mentioned in the Introduction, there is a relationship between
`characterizations of regulated functions' and `pointwise selection principles'.
The second main result is a pointwise selection principle for
uniform space valued functions in terms of the approximate variation:

\begin{thm} \label{th:2}
Suppose $\varnothing\ne T\subset\RB$ and $(X,\UC)$ is a Hausdorff uniform space
with an at most countable gage of pseudometrics $\{d_p\}_{p\in\PC}$ of the
uniformity~$\UC$. If\/ $\{f_j\}\subset X^T$ is a pointwise relatively sequentially compact
sequence of functions on $T$ such that
  \begin{equation} \label{e:Vepfj}
\limsup_{j\to\infty}V_{\vep,\,p}(f_j,T)<\infty\,\,\,\mbox{for all}\,\,\,\vep>0
\mbox{ and } p\in\PC,
  \end{equation}
then there is a subsequence of $\{f_j\}$ which converges pointwise on $T$ to a
function $f\in\mbox{\rm B}(T;X)\cap\mbox{\rm Reg}(T;X)$.
\end{thm}

The novelty of this theorem is threefold. First, functions $f_j$ take their values in a
uniform space, which (as it was seen earlier) is more general than a metric space.
Second, condition \eq{e:Vepfj} is {\em necessary\/} for uniformly convergent
sequences $\{f_j\}$ (but not for pointwise convergent sequences $\{f_j\}$). Third,
Theorem~\ref{th:2} may be applied to sequences $\{f_j\}$ of non-regulated functions.
These issues are addressed in Section~\ref{s:proofs}.

\section{Properties of the approximate variation} \label{s:prappvar}
 
\subsection{Bounded functions and functions of bounded variation.}

Given $f,g\in X^T$, $s,t\in T$, and $p\in\PC$, by the triangle inequality for $d_p$,
we find
  \begin{equation} \label{e:2_1}
d_{T,\,p}(f,g)\le|f(T)|_p+d_p(f(t),g(t))+|g(T)|_p
  \end{equation}
and
  \begin{equation} \label{e:2-2}
d_p(f(s),f(t))\le d_p(g(s),g(t))+2d_{T,\,p}(f,g)
  \end{equation}
and so, the definition of the $d_p$-oscillation and \eq{e:2-2} imply
  \begin{equation} \label{e:2^3}
|f(T)|_p\le|g(T)|_p+2d_{T,\,p}(f,g).
  \end{equation}
By \eq{e:2_1} and \eq{e:2^3}, we have $d_{T,\,p}(f,g)<\infty$ for all
$f,g\in\mbox{\rm B}_p(T;X)$ with $p\in\PC$ and, given a constant function $c\in X^T$,
\mbox{$\mbox{\rm B}_p(T;X)=\{f\in X^T:d_{T,\,p}(f,c)<\infty\}$} for all $p\in\PC$,
and so,
$\mbox{\rm B}(T;X)=\{f\in X^T:\mbox{$d_{T,\,p}(f,c)<\infty$ for all $p\in\PC$}\}$.

\smallbreak
Clearly,
  \begin{equation} \label{e:leVp}
|f(T)|_p\le V_p(f,T)\!\quad\mbox{for all}\!\quad\!f\in X^T\!\quad\!
\mbox{and}\!\quad\! p\in\PC,
  \end{equation}
and so,
  \begin{equation*}
\mbox{\rm BV}_p(T;X)\subset\mbox{\rm B}_p(T;X)\quad\mbox{and}\quad
\mbox{\rm BV}(T;X)\subset\mbox{\rm B}(T;X).
  \end{equation*}

Given $f\in X^T$ and $p\in\PC$, the functional $V_p(\cdot,\cdot)$ has the following
two properties: (i) {\em additivity\/} of $V_p(f,\,\cdot\,)$ in the second variable
(cf. \cite{MatSb}, \cite{JMAA19}):
  \begin{equation} \label{e:Vpadd}
V_p(f,T)=V_p(f,T\cap(-\infty,t])+V_p(f,T\cap[t,\infty))\quad\mbox{for \,all \,\,$t\in T$;}
  \end{equation}
(ii) sequential {\em lower semicontinuity\/} of $V_p(\,\cdot\,,T)$ in the first variable
 (cf. \cite{MatTr}):
  \begin{equation} \label{e:Vplsc}
\mbox{if $\{f_j\}\subset X^T$ and $f_j\to_p f$ on $T$, then
$V_p(f,T)\le\liminf\limits_{j\to\infty}V_p(f_j,T)$.}
  \end{equation}

\subsection{The approximate variation.}

The essential properties of the approximate variation are gathered in

\begin{lem} \label{l:ess}
Given $f\in X^T$ and $p\in\PC$, we have\/{\rm:}\par\vspace{-8pt}
  \begin{itemize}
\renewcommand{\itemsep}{-1.0pt plus 0.5pt minus 0.25pt}
\item[{\rm(a)}] the function $\vep\mapsto V_{\vep,\,p}(f,T):(0,\infty)\to[0,\infty]$ is
  nonincreasing, and so,
  $V_{\vep+0,\,p}(f,T)\le V_{\vep,\,p}(f,T)\le V_{\vep-0,\,p}(f,T)$
  for all $\vep>0;$%
\footnote{As usual, $V_{\vep+0,\,p}(f,T)$ and $V_{\vep-0,\,p}(f,T)$ are the limits
from the right and from the left of the function $\vep\mapsto V_{\vep,\,p}(f,T)$, respectively.}
\item[{\rm(b)}] if\/ $\varnothing\ne T_1\subset T_2\subset T$ and $\vep>0$, then
  $V_{\vep,\,p}(f,T_1)\le V_{\vep,\,p}(f,T_2);$
\item[{\rm(c)}] $\lim\limits_{\vep\to+0}V_{\vep,\,p}(f,T)=%
  \sup\limits_{\vep>0}V_{\vep,\,p}(f,T)=V_p(f,T);$
\item[{\rm(d)}] $|f(T)|_p\le V_{\vep,\,p}(f,T)+2\vep$ for all $\vep>0;$
\item[{\rm(e)}] if\/ $\vep>0$, $t\in T$, $T_t^-=T\cap(-\infty,t]$ and\/
  $T_t^+=T\cap[t,\infty)$, then
  \begin{equation*}
V_{\vep,\,p}(f,T_t^-)+V_{\vep,\,p}(f,T_t^+)\le V_{\vep,\,p}(f,T)\le
V_{\vep,\,p}(f,T_t^-)+V_{\vep,\,p}(f,T_t^+)+2\vep.
  \end{equation*} 
  \end{itemize}
\end{lem}

\begin{pf}
(a) Suppose $0<\vep_1<\vep_2$. If $g\in X^T$ and $d_{T,\,p}(f,g)\le\vep_1$, then
$d_{T,\,p}(f,g)\le\vep_2$, and so, definition \eq{d:Vep} implies
$V_{\vep_2,\,p}(f,T)\le V_{\vep_1,\,p}(f,T)$.

\smallbreak
(b) Since $T_1\subset T_2$, we have $d_{T_1,\,p}(f,g)\le d_{T_2,\,p}(f,g)$ for all
$g\in X^T$, whence, by definition \eq{d:Vep},
$V_{\vep,\,p}(f,T_1)\le V_{\vep,\,p}(f,T_2)$.

\smallbreak
(c) By (a), the quantity $C_p=\lim_{\vep\to+0}V_{\vep,\,p}(f,T)$, which actually
 is equal to $\sup_{\vep>0}V_{\vep,\,p}(f,T)$, is well-defined in $[0,\infty]$. Assume first
that $f\in\mbox{\rm BV}_p(T;X)$. By \eq{d:Vep},
$V_{\vep,\,p}(f,T)\le V_p(f,T)$ for all $\vep>0$, and so, $C_p\le V_p(f,T)<\infty$.
To prove the reverse inequality, we apply the definition of $C_p$: for every $\eta>0$
there is $\delta=\delta(\eta)>0$ such that $V_{\vep,\,p}(f,T)\le C_p+\eta$ for all
$\vep\in(0,\delta)$. If $\{\vep_k\}_{k=1}^\infty$ is a sequence in $(0,\delta)$ such that
$\vep_k\to0$ as $k\to\infty$, then, for every $k\in\NB$, the definition of
$V_{\vep_k,\,p}(f,T)$ yields the existence of $g_k\in\mbox{\rm BV}_p(T;X)$ such that
$d_{T,\,p}(f,g_k)\le\vep_k$ and $V_p(g_k,T)\le C_p+\eta$. Since $\vep_k\to0$,
$g_k\rra_p f$ on~$T$, hence $g_k\to_p f$ on $T$, and so, by \eq{e:Vplsc},
  \begin{equation*}
V_p(f,T)\le\liminf_{k\to\infty}V_p(g_k,T)\le C_p+\eta.
  \end{equation*}
By the arbitrariness of $\eta>0$, $V_p(f,T)\le C_p$. Thus, $C_p$ and $V_p(f,T)$ are
finite or not finite simultaneously, and $C_p=V_p(f,T)$, which proves~(c).

\smallbreak
(d) We may assume that $\vep>0$ is such that $V_{\vep,\,p}(f,T)$ is finite.
By definition \eq{d:Vep}, given $\eta>0$, there is $g=g_\eta\in\mbox{\rm BV}_p(T;X)$
such that $d_{T,\,p}(f,g)\le\vep$ and $V_p(g,T)\le V_{\vep,\,p}(f,T)+\eta$. Now,
\eq{e:2^3} and \eq{e:leVp} imply
  \begin{equation*}
|f(T)|_p\le|g(T)|_p+2d_{T,\,p}(f,g)\le V_p(g,T)+2\vep\le V_{\vep,\,p}(f,T)+\eta+2\vep,
  \end{equation*}
and the inequality in (d) follows due to the arbitrariness of $\eta>0$.

\smallbreak
(e) First, we prove the left-hand side inequality, in which we may assume that
$V_{\vep,\,p}(f,T)$ is finite. By definition \eq{d:Vep}, for every $\eta>0$ there is
$g=g_\eta\in\mbox{\rm BV}_p(T;X)$ such that $d_{T,\,p}(f,g)\le\vep$ and
$V_p(g,T)\le V_{\vep,\,p}(f,T)+\eta$. We set $g^-(s)=g(s)$ for $s\in T_t^-$ and
$g^+(s)=g(s)$ for $s\in T_t^+$, and note that $g^-(t)=g(t)=g^+(t)$. Since
$g^-\in\mbox{\rm BV}_p(T_t^-;X)$, $g^+\in\mbox{\rm BV}_p(T_t^+;X)$,
  $$d_{T_t^-\!,\,p}(f,g^-)\le d_{T,\,p}(f,g)\le\vep\quad\mbox{and}\quad
     d_{T_t^+\!,\,p}(f,g^+)\le d_{T,\,p}(f,g)\le\vep,$$
definition \eq{d:Vep} implies $V_{\vep,\,p}(f,T_t^-)\le V_p(g^-,T_t^-)$ and
$V_{\vep,\,p}(f,T_t^+)\le V_p(g^+,T_t^+)$, and so, by the additivity property
\eq{e:Vpadd}, we find
  \begin{align*}
V_{\vep,\,p}(f,T_t^-)+V_{\vep,\,p}(f,T_t^+)&\le V_p(g^-,T_t^-)+V_p(g^+,T_t^+)
  =V_p(g,T_t^-)+V_p(g,T_t^+)\\[3pt]
&=V_p(g,T)\le V_{\vep,\,p}(f,T)+\eta\quad\mbox{for all \,$\eta>0$.}
  \end{align*}
This establishes the left-hand side inequality.

\smallbreak
In order to prove the right-hand side inequality, we may assume that
$V_{\vep,\,p}(f,T_t^-)$ and $V_{\vep,\,p}(f,T_t^+)$ are finite.
Note that if $T_t^-\!=\!\{t\}$, then $T_t^+\!=\!T$, and if $T_t^+\!=\!\{t\}$,
 then $T_t^-\!=\!T$;
in both these cases the inequality is clear. Assume that $T_t^-\ne\{t\}$ and
$T_t^+\ne\{t\}$, so that $T\cap(-\infty,t)\!\ne\!\varnothing$ and
$T\cap(t,\infty)\!\ne\!\varnothing$. By definition \eq{d:Vep}, for every ${\eta^-}>0$ and
${\eta^+}>0$ there are $g^-\in\mbox{\rm BV}_p(T_t^-;X)$ and
$g^+\in\mbox{\rm BV}_p(T_t^+;X)$ with the properties: 
$d_{T_t^-\!,\,p}(f,g^-)\!\le\!\vep$, \mbox{$d_{T_t^+\!,\,p}(f,g^+)\!\le\!\vep$},
$V_p(g^-,T_t^-)\le V_{\vep,\,p}(f,T_t^-)+{\eta^-}$ and
$V_p(g^+,T_t^+)\le V_{\vep,\,p}(f,T_t^+)+{\eta^+}$. Given $x\in X$ (to be
specified below), we define $g\in\mbox{\rm BV}_p(T;X)$ by
  $$g(s)=\left\{\begin{array}{ccl}
    \!\!g^-(s) & \!\!\mbox{if}\! & s\in T\cap(-\infty,t),\\
    \!\!x & \!\!\mbox{if}\! & s=t,\\
    \!\!g^+(s) & \!\!\mbox{if}\! & s\in T\cap(t,\infty).
    \end{array}\right.$$
Suppose that the following two inequalities are already established:
  \begin{equation} \label{e:gT-}
V_p(g,T_t^-)\le V_p(g^-,T_t^-)+d_p(g(t),g^-(t))
  \end{equation}
and
  \begin{equation} \label{e:gT+}
V_p(g,T_t^+)\le V_p(g^+,T_t^+)+d_p(g(t),g^+(t)).
  \end{equation}
By the additivity property \eq{e:Vpadd} of $V_p(\cdot,\cdot)$, these inequalities imply
  \begin{align}
V_p(g,T)&=V_p(g,T_t^-)+V_p(g,T_t^+) \nonumber \\[3pt]
&\le V_{\vep,\,p}(f,T_t^-)+{\eta^-}+d_p(x,g^-(t)) \nonumber \\[3pt]
&\quad\,+V_{\vep,\,p}(f,T_t^+)+{\eta^+}+d_p(x,g^+(t)).\label{e:x+g}
  \end{align}
Now, we put $x=g^-(t)$ (by symmetry, we may have put $x=g^+(t)$ as well).
Since $g=g^-$ on $T_t^-$ and $g=g^+$ on $T\cap(t,\infty)\subset T_t^+$, we get
  \begin{equation} \label{e:max2}
d_{T,\,p}(f,g)\le\max\bigl\{d_{T_t^-\!,\,p}(f,g^-),d_{T_t^+\!,\,p}(f,g^+)\bigr\}\le\vep.
  \end{equation}
Noting that (cf. \eq{e:x+g})
  \begin{align*}
d_p(x,g^+(t))&=d_p(g^-(t),g^+(t))\le d_p(g^-(t),f(t))+d_p(f(t),g^+(t))\\[3pt]
&\le d_{T_t^-\!,\,p}(g^-,f)+d_{T_t^+\!,\,p}(f,g^+)\le\vep+\vep=2\vep,
  \end{align*}
we conclude from \eq{d:Vep}, \eq{e:x+g} and \eq{e:max2} that
  \begin{equation*}
V_{\vep,\,p}(f,T)\le V_p(g,T)\le V_{\vep,\,p}(f,T_t^-)+{\eta^-}+V_{\vep,\,p}(f,T_t^+)
+{\eta^+}+2\vep,
  \end{equation*}
and the inequality in (e) follows due to the arbitrariness of ${\eta^-}\!>\!0$
and ${\eta^+}\!>\!0$.

\smallbreak
It remains to establish \eq{e:gT-} (since \eq{e:gT+} is similar).
Let $\{t_i\}_{i=0}^m\subset T_t^-$ be a partition of $T_t^-$, i.e.,
$t_0\le t_1\le\ldots\le t_{m-1}<t_m=t$. Since $g(s)=g^-(s)$ for all $s\in T$ with $s<t$,
we have
  \begin{align*}
\sum_{i=1}^m\!d_p(g(t_i),g(t_{i-1}))&=\sum_{i=1}^{m-1}\!d_p(g(t_i),g(t_{i-1}))
  +d_p(g(t_m),g(t_{m-1}))\\
&=\sum_{i=1}^{m-1}\!d_p(g^-(t_i),g^-(t_{i-1}))+d_p(g^-(t_m),g^-(t_{m-1}))\\
&\qquad+d_p(g(t_m),g(t_{m-1}))\!-\!d_p(g^-(t_m),g^-(t_{m-1}))\\[5pt]
&\le V_p(g^-,T_t^-)\!+\!\bigl|d_p(g(t),g^-(t_{m-1}))\!-\!d_p(g^-(t),g^-(t_{m-1}))
  \bigr|\\[5pt]
&\le V_p(g^-,T_t^-)+d_p(g(t),g^-(t)),
  \end{align*}
where the last inequality is due to the triangle inequality for $d_p$. Taking the supremum
over all partitions of $T_t^-$, we obtain inequality \eq{e:gT-}.
\qed\end{pf}

\begin{rem}
Inequalities in Lemma~\ref{l:ess}\,(e) are sharp (cf.\ \cite[Example~3.2]{arXiv19}).
\end{rem}

\subsection{Regulated functions and Dirichlet's function.}

Now, we are in a position to prove Theorem~\ref{th:1}.

\begin{pthone}
Suppose $p\in\PC$, and $f\in X^T$ is such that $V_{\vep,\,p}(f,T)$ is finite for all $\vep>0$.
Given $\tau\in T$, which is a left limit point for $T$, let us show that $d_p(f(s),f(t))\to0$
as $T\ni s,t\to\tau-0$ (the arguments for $\tau'\in T$ being a right limit point for $T$
are similar). For any $\vep>0$, we define the {\em$\vep$-$d_p$-variation function\/} by
$\vfi_{\vep,\,p}(t)=V_{\vep,\,p}(f,T_t^-)$, $t\in T$, where $T_t^-=T\cap(-\infty,t]$.
By Lemma~\ref{l:ess}\,(b), given $s,t\in T$ with $s\le t$ (and so, $T_s^-\subset T_t^-$),
we find
  \begin{equation*}
0\le\vfi_{\vep,\,p}(s)\le\vfi_{\vep,\,p}(t)\le V_{\vep,\,p}(f,T)<\infty,
  \end{equation*}
i.e., $\vfi_{\vep,\,p}:T\to[0,\infty)$ is a bounded and nondecreasing function.
Consequently, the left limit
  \begin{equation*}
\mbox{$\lim\limits_{T\ni t\to\tau-0}\vfi_{\vep,\,p}(t)=\sup\limits_{t\in T\cap(-\infty,\tau)}
\vfi_{\vep,\,p}(t)$ \,exists \,in \,$[0,\infty)$.}
  \end{equation*}
It follows that there is $\delta=\delta(\vep)>0$ such that
$|\vfi_{\vep,\,p}(t)-\vfi_{\vep,\,p}(s)|<\vep$ for all $s,t\in T\cap(\tau-\delta,\tau)$.
Now, assume that $s,t\in T\cap(\tau-\delta,\tau)$ are arbitrary such that $s<t$.
Lemma~\ref{l:ess}\,(e) (in which $T$ is replaced by $T_t^-$, so that
$(T_t^-)_s^-=T_s^-$ and $(T_t^-)_s^+=T\cap[s,t]$ for $s\in T_t^-$) implies
  \begin{equation*}
V_{\vep,\,p}(f,T\cap[s,t])\le V_{\vep,\,p}(f,T_t^-)-V_{\vep,\,p}(f,T_s^-)
=\vfi_{\vep,\,p}(t)-\vfi_{\vep,\,p}(s)<\vep.
  \end{equation*}
By the definition of $V_{\vep,\,p}(f,T\cap[s,t])$, there is
$g=g_{\vep,\,p}\in\mbox{\rm BV}_p(T\cap[s,t];X)$ such that
$d_{T\cap[s,t],\,p}(f,g)\le\vep$ and $V_p(g,T\cap[s,t])\le\vep$, and so,
\eq{e:2-2} yields
  \begin{equation*}
d_p(f(s),f(t))\le d_p(g(s),g(t))+2d_{T\cap[s,t],\,p}(f,g)\le V_p(g,T\cap[s,t])+2\vep
\le3\vep.
  \end{equation*}
This completes the proof of the equality
$\lim_{T\ni s,t\to\tau-0}d_p(f(s),f(t))=0$.

\smallbreak
Suppose now that $f\in\mbox{\rm Reg}_p(I;X)$ with $I=[a,b]$. By virtue of
\cite[Lemma~4]{MatTr}, there is a sequence of step functions%
\footnote{Recall that $g\in X^I$ is a {\em step function\/} if, for some partition
$a=t_0<t_1<t_2<\ldots<t_{m-1}<t_m=b$ of $I=[a,b]$, $g$ takes a constant value
on each open interval $(t_{i-1},t_i)$, $i=1,2,\ldots,m$. Clearly,
$g\in\mbox{\rm BV}_p(I;X)$ for all $p\in\PC$.}
$\{g_j\}\subset X^I$ such that $g_j\rra_p f$ on $I$, and so, given $\vep>0$,
there is $j(\vep)\in\NB$ such that $d_{I,\,p}(f,g_{j(\vep)})\le\vep$. Since
$g_{j(\vep)}\in\mbox{\rm BV}_p(I;X)$, definition \eq{d:Vep} yields
$V_{\vep,\,p}(f,I)\le V_p(g_{j(\vep)},I)<\infty$, which was to be proved.
\qed\end{pthone}

As an illustration of Theorem~\ref{th:1}, we present an example.

\begin{exa} \label{ex:Dir}
Let $I=[a,b]$, $x,y\in X$, $x\ne y$, and $\QB$ be the set of all rational numbers.
We denote by $f=\DC_{x,y}\in X^I$ the {\em Dirichlet function\/} given~by
  \begin{equation} \label{e:Dir}
\mbox{$\DC_{x,y}(t)=x$ \,if \,$t\in I\cap\QB$, \,and \,$\DC_{x,y}(t)=y$ \,if
\,$t\in I\setminus\QB$.}
  \end{equation}

\smallbreak
It is clearly seen that $f\notin\mbox{\rm Reg}(I;M)$: in fact, since $X$ is Hausdorff
and $x\ne y$, there is $p_0\in\PC$ such that $d_{p_{{}_0}}(x,y)>0$, and so, if, say,
$a<\tau\le b$, then for all $\delta\in(0,\tau-a)$, $s\in(\tau-\delta,\tau)\cap\QB$
and $t\in(\tau-\delta,\tau)\setminus\QB$, we have
$d_{p_{{}_0}}(f(s),f(t))=d_{p_{{}_0}}(x,y)>0$,
so that the limit of $d_{p_{{}_0}}(f(s),f(t))$ as $I\ni s,t\to\tau-0$ does not exist.

\smallbreak
We claim that, for $\vep>0$ and all $p\in\PC$ with $d_p(x,y)>0$, we have
  \begin{equation} \label{e:inf0}
V_{\vep,\,p}(f,I)=\left\{\begin{array}{ccl}
    \!\!\infty & \!\!\mbox{if}\! & \vep<d_p(x,y)/2,\\[4pt]
    \!\!0 & \!\!\mbox{if}\! & \vep\ge d_p(x,y)\\
    \end{array}\right.
  \end{equation}
(the value $V_{\vep,\,p}(f,I)$ for $d_p(x,y)/2\le\vep<d_p(x,y)$ depends on the
structure of $X$ such as generalized convexity, cf. \cite[Example~1]{Studia}).
In fact, if \mbox{$\vep\ge d_p(x,y)$}, we may set $g(t)\!=\!x$ (or $g(t)\!=\!y$)
for all $t\!\in\! I$, so that $d_{I,\,p}(f,g)\!=\!d_p(x,y)\!\le\!\vep$ and, by \eq{d:Vep},
$0\le V_{\vep,\,p}(f,I)\le V_p(g,I)=0$. Suppose $0<\vep<d_p(x,y)/2$.
We are going to show that $V_p(g,I)=\infty$ for all $g\in X^I$ with $d_{I,\,p}(f,g)\le\vep$,
which, in accordance with \eq{d:Vep}, will imply $V_{\vep,\,p}(f,I)=\inf\varnothing=\infty$.
Let $m\in\NB$ and $a\le t_0<t_1<t_2<\ldots<t_{m-1}<t_m\le b$ be a partition of
$I=[a,b]$ such that $\{t_{2i}\}_{i=0}^m\subset I\cap\QB$ and
$\{t_{2i-1}\}_{i=1}^m\subset I\setminus\QB$. Given $i\in\{1,2,\dots,m\}$,
\eq{e:2-2} implies
  $$d_p(x,y)=d_p(f(t_{2i}),f(t_{2i-1}))\le d_p(g(t_{2i}),g(t_{2i-1}))+2\vep,$$
whence, by the definition of $V_p(g,I)$,
  $$V_p(g,I)\ge\sum_{i=1}^m d_p(g(t_{2i}),g(t_{2i-1}))\ge(d_p(x,y)-2\vep)m\quad
     \,\forall\,m\in\NB.$$
\end{exa}

\subsection{Uniform convergence vs.\ pointwise convergence.}

The next lemma shows the interplay of the approximate variation and the uniform
convergence.

\begin{lem} \label{l:unif}
Let $\varnothing\ne T\subset\RB$ and $(X,\UC)$ be a Hausdorff uniform space.
Suppose $p\in\PC$, $f\in X^T$, $\{f_j\}\subset X^T$, and $f_j\rra_p f$ on~$T$.
We have\/{\rm:} \par\vspace{-2pt}
  \begin{itemize}
\renewcommand{\itemsep}{-1.0pt plus 0.5pt minus 0.25pt}
\item[{\rm(a)}] $\!\displaystyle V_{\vep+0,\,p}(f,T)\!\le\!
  \liminf_{j\to\infty}V_{\vep,\,p}(f_j,T)\!\le\!\limsup_{j\to\infty} V_{\vep,\,p}(f_j,T)
  \!\le\! V_{\vep-0,\,p}(f,T)$, \mbox{$\vep\!>\!0;$}
\item[{\rm(b)}] $\!$if\/ $V_{\vep,\,p}(f_j,T)\!<\!\infty$ for all $\vep\!>\!0$ and
  $j\!\in\!\NB$, then $V_{\vep,\,p}(f,T)\!<\!\infty$~for~all~\mbox{$\vep\!>\!0$}.
  \end{itemize}
\end{lem}

\begin{pf}
(a) Only the first and the last inequalities are to be established.

\smallbreak
1. In order to prove the first inequality, we may assume (passing to a suitable subsequence
of $\{f_j\}$ if necessary) that its right-hand side is equal to
$C_p=\lim_{j\to\infty}V_{\vep,\,p}(f_j,T)<\infty$. Let $\eta>0$ be given arbitrarily.
Then, there is $j_0=j_0(\eta)\in\NB$ such that $V_{\vep,\,p}(f_j,T)\le C_p+\eta$
for all $j\ge j_0$. By the definition of $V_{\vep,\,p}(f_j,T)$, for every $j\ge j_0$ there is
$g_j\in\mbox{\rm BV}_p(T;X)$ (also depending on $\eta$ and $p$) such that
$d_{T,\,p}(f_j,g_j)\!\le\!\vep$ and $V_p(g_j,T)\!\le\! V_{\vep,\,p}(f_j,T)\!+\!\eta$. Since
\mbox{$f_j\rra_p f$} on $T$, we find $d_{T,\,p}(f_j,f)\to0$ as $j\to\infty$, and so, there is
$j_1=j_1(\eta)\in\NB$ such that $d_{T,\,p}(f_j,f)\le\eta$ for all $j\ge j_1$.
Taking into account that
  $$d_{T,\,p}(f,g_j)\le d_{T,\,p}(f,f_j)+d_{T,\,p}(f_j,g_j)\le\eta+\vep\quad
     \mbox{for all \,$j\ge\max\{j_0,j_1\}$,}$$
and applying definition \eq{d:Vep}, we get
  $$V_{\eta+\vep,\,p}(f,T)\le V_p(g_j,T)\le V_{\vep,\,p}(f_j,T)+\eta\le
     (C_p+\eta)+\eta=C_p+2\eta.$$
Passing to the limit as $\eta\to+0$, we arrive at $V_{\vep+0,\,p}(f,T)\le C_p$, which
was to be proved.

\smallbreak
2. To establish the third inequality in (a), we may assume (with no loss of generality) that
$V_{\vep-0,\,p}(f,T)$ is finite. Given $\eta>0$, there is
\mbox{$\delta=\delta(\eta,\vep)\in(0,\vep)$}
such that $V_{\vep'\!,\,p}(f,T)\le V_{\vep-0,\,p}(f,T)+\eta$ for all
$\vep'\in(\vep-\delta,\vep)$. Since $f_j\rra_p f$ on $T$, for every
$\vep'\in(\vep-\delta,\vep)$ there is $j_0=j_0(\vep',\vep)\in\NB$ such that
$d_{T,\,p}(f_j,f)\le\vep-\vep'$ for all $j\ge j_0$. By the definition of $V_{\vep'\!,\,p}(f,T)$,
given $j\in\NB$, there is $g_j\in\mbox{\rm BV}_p(T;X)$ (also depending on $\vep'$
and $p$) such that $d_{T,\,p}(f,g_j)\le\vep'$ and
  $$V_{\vep'\!,\,p}(f,T)\le V_p(g_j,T)\le V_{\vep'\!,\,p}(f,T)+(1/j).$$
Hence, $\lim_{j\to\infty}V_p(g_j,T)=V_{\vep'\!,\,p}(f,T)$. Noting that, for all $j\ge j_0$,
  $$d_{T,\,p}(f_j,g_j)\le d_{T,\,p}(f_j,f)+d_{T,\,p}(f,g_j)\le(\vep-\vep')+\vep'=\vep,$$
we find from \eq{d:Vep} that $V_{\vep,\,p}(f_j,T)\le V_p(g_j,T)$ for all $j\ge j_0$. Thus,
  $$\limsup_{j\to\infty}V_{\vep,\,p}(f_j,T)\le\lim_{j\to\infty}V_p(g_j,T)=
     V_{\vep'\!,\,p}(f,T)\le V_{\vep-0,\,p}(f,T)+\eta.$$
It remains to take into account the arbitrariness of $\eta>0$.

\smallbreak
(b) Let $\vep>0$ and $0<\vep'<\vep$. Given $j\in\NB$, since
$V_{\vep'\!,\,p}(f_j,T)<\infty$, definition \eq{d:Vep} implies the existence of
$g_j\in\mbox{\rm BV}_p(T;X)$ such that $d_{T,\,p}(f_j,g_j)\le\vep'$ and
$V_p(g_j,T,)\le V_{\vep'\!,\,p}(f_j,T)+1$. Since $f_j\rra_p f$ on $T$, there is
\mbox{$j_0=j_0(\vep,\vep')\in\NB$} such that $d_{T,\,p}(f_{j_{{}_0}},f)\le\vep-\vep'$.
Now, it follows from
  $$d_{T,\,p}(f,g_{j_{{}_0}})\le d_{T,\,p}(f,f_{j_{{}_0}})+
     d_{T,\,p}(f_{j_{{}_0}},g_{j_{{}_0}})\le(\vep-\vep')+\vep'=\vep$$
and definition \eq{d:Vep} that $V_{\vep,\,p}(f,T)\le V_p(g_{j_{{}_0}},T)\le%
V_{\vep'\!,\,p}(f_{j_{{}_0}},T)+1<\infty$.
\qed\end{pf}

A few comments and examples concerning Lemma~\ref{l:unif} are in order.

\begin{exa} \label{ex:vez}
The left limit $V_{\vep-0,\,p}(f,T)$ in Lemma~\ref{l:unif}\,(a) cannot in general be
replaced by $V_{\vep,\,p}(f,T)$. We demonstrate this in the case of a normed linear
space $(X,\|\cdot\|)$ equipped with the canonical metric $d(x,y)=\|x-y\|$, $x,y\in X$
(here we omit the subscript $p$ in the notations $d_T=d_{T,\,p}$, $V=V_p$ and
$V_\vep=V_{\vep,\,p}$). Let $T=I=[a,b]$ and $\{x_j\}$, $\{y_j\}\subset X$ be two
sequences such that $\|x_j-x\|\to0$ and $\|y_j-y\|\to0$ as $j\to\infty$, where
$x,y\in X$, $x\ne y$. If $f_j=\DC_{x_j,y_j}$, $j\in\NB$, and $f=\DC_{x,y}$ are
Dirichlet functions \eq{e:Dir} on $I$, then $f_j\rra f$ on $I$, because
  $$d_I(f_j,f)\!=\!\sup_{t\in I}\|f_j(t)\!-\!f(t)\|\!=\!\max\bigl\{\|x_j\!-\!x\|,\|y_j\!-\!y\|\bigr\}
    \to0\quad\mbox{as}\quad j\to\infty.$$
According to \eq{e:inf0}, the values $V_\vep(f,I)$ are given for $\vep>0$ by
  $$\mbox{$V_\vep(f,I)=\infty$ if $\vep<\|x\!-\!y\|/2$, \,and $V_\vep(f,I)=0$ if
     $\vep\ge\|x\!-\!y\|/2$}$$
(note that if $\|x-y\|/2\le\vep<\|x-y\|$, we may set $g(t)=(x+y)/2$ for all $t\in I$,
which implies $d_I(f,g)=\|x-y\|/2\le\vep$, and so, $0\le V_\vep(f,I)\le V(g,I)=0$).
Similarly, for (large) $j\in\NB$,
  $$\mbox{$V_\vep(f_j,I)=\infty$ if $\vep<\|x_j\!-\!y_j\|/2$, \,and $V_\vep(f_j,I)=0$ if
     $\vep\ge\|x_j\!-\!y_j\|/2$.}$$
Setting $\vep=\|x-y\|/2$, we find
  $$V_{\vep+0}(f,I)=V_\vep(f,I)=0<\infty=V_{\vep-0}(f,I),$$
whereas, if $x_j=\alpha_jx$ and $y_j=\alpha_jy$ with $\alpha_j=1+(1/j)$, $j\in\NB$, we get
$\vep<\alpha_j\|x-y\|/2=\|x_j-y_j\|/2$, which implies $V_\vep(f_j,I)=\infty$ for all
$j\in\NB$, and so, $\lim_{j\to\infty}V_\vep(f_j,I)=\infty$.
\end{exa}

\begin{exa}
In this example, we show that the right-hand side inequality in Lemma~\ref{l:unif}\,(a)
does not hold in general if $\{f_j\}\subset X^T$ converges to $f\in X^T$ only pointwise
on~$T$. Suppose $p_0\in\PC$ is such that
$C_{p_{{}_0}}=\inf_{j\in\NB}|f_j(T)|_{p_{{}_0}}>0$  and $f=c$ (is a constant function)
on~$T$. Given $0<\vep<C_{p_{{}_0}}/2$ and $j\in\NB$, Lemma~\ref{l:ess}\,(d) implies
  $$V_{\vep,\,p_{{}_0}}(f_j,T)\ge|f_j(T)|_{p_{{}_0}}-2\vep\ge C_{p_{{}_0}}-2\vep
     >0=V_{\vep,p_{{}_0}}(c,T)=V_{\vep-0,p_{{}_0}}(f,T).$$

\smallbreak
More specifically, let $x,y\in X$, $x\ne y$, $p_0\in\PC$ be such that $d_{p_{{}_0}}(x,y)>0$,
and $\{\tau_j\}_{j=1}^\infty\subset(a,b)\subset T=[a,b]$ be such that $\tau_j\to a$
as $j\to\infty$. Defining $f_j\in X^T$ by $f_j(\tau_j)=x$ and $f_j(t)=y$ if
$t\in T\setminus\{\tau_j\}$, we get $C_{p_{{}_0}}=d_{p_{{}_0}}(x,y)>0$ and
$f_j\to_p c\equiv y$ on $T$ for all $p\in\PC$.

\smallbreak
Note that the arguments above are not valid for the uniform convergence: in fact,
if $p\in\PC$ and $f_j\rra_p f=c$ on $T$, then, by \eq{e:2^3},
  $$|f_j(T)|_p\le2d_{T,\,p}(f_j,c)\to0\quad\mbox{as}\quad j\to\infty,$$
and so, $C_p\equiv\inf_{j\in\NB}|f_j(T)|_p=0$.
\end{exa}

\begin{exa} \label{exa4}
Lemma~\ref{l:unif}\,(b) is wrong for the pointwise convergence of $\{f_j\}$ to~$f$.
To see this, let $T=I=[a,b]$, $x,y\in X$, $x\ne y$, and, given $j\in\NB$, define
$f_j\in X^I$ at $t\in I$ by: $f_j(t)=x$ if $j!\cdot t$ is integer, and $f_j(t)=y$ otherwise.
Each $f_j$ is a step function on $I$, so $f_j\in\mbox{\rm Reg}(I;X)$ and, hence,
by Theorem~\ref{th:1}, $V_{\vep,\,p}(f_j,I)<\infty$ for all $\vep>0$ and $p\in\PC$.
Clearly, $f_j$ converges (only) pointwise on $I$ to the Dirichlet function $f=\DC_{x,y}$
from \eq{e:Dir}, and so, by \eq{e:inf0}, $V_{\vep,\,p}(f,I)=\infty$ for all $p\in\PC$
with $d_p(x,y)>0$ and $0<\vep<d_p(x,y)/2$.
\end{exa}

\section{Proof of the main result} \label{s:proofs}

We denote by $\mbox{\rm Mon}(T;\RB^+)$ the family of all bounded nondecreasing
functions mapping $T$ into $\RB^+=[0,\infty)$ (where $\RB^+$ may be replaced
by~$\RB$). It is worthwhile to recall the classical Helly selection principle for an arbitrary
set $T\subset\RB$ (cf.\ \cite[Proof of Theorem~1.3]{Sovae}): {\em a uniformly
bounded sequence of functions from $\mbox{\rm Mon}(T;\RB)$ contains a subsequence
which converges pointwise on $T$ to a function from $\mbox{\rm Mon}(T;\RB)$}.

\begin{pthscnd}
By Lemma~\ref{l:ess}\,(b), given $\vep>0$, $p\in\PC$ and $j\in\NB$, the
$\vep$-$d_p$-variation function $t\mapsto V_{\vep,\,p}(f_j,T_t^-)$ is nondecreasing
on~$T$ (recall that $T_t^-=T\cap(-\infty,t]$ for $t\in T$). By assumption \eq{e:Vepfj},
for every $\vep>0$ and $p\in\PC$ there are $j^*(\vep,p)\in\NB$ and $C(\vep,p)>0$
such that $V_{\vep,\,p}(f_j,T)\le C(\vep,p)$ for all natural $j\ge j^*(\vep,p)$.
Again by Lemma~\ref{l:ess}\,(b),
  $$V_{\vep,\,p}(f_j,T_t^-)\le V_{\vep,\,p}(f_j,T)\le C(\vep,p)\quad
     \mbox{for all $t\in T$ and $j\ge j^*(\vep,p)$,}$$
and so, the sequence of nondecreasing functions
$\{t\mapsto V_{\vep,\,p}(f_j,T_t^-)\}_{j=j^*(\vep,p)}^\infty$ is uniformly
bounded in $t\in T$ by constant $C(\vep,p)$.

\smallbreak
We divide the rest of the proof into four steps. With no loss of generality we assume
that $\PC=\NB$.

\smallbreak
1. Making use of the Cantor diagonal process, let us show that for each decreasing
sequence $\{\vep_k\}_{k=1}^\infty$ of positive numbers $\vep_k\to0$ (as $k\to\infty$)
there is a subsequence of $\{f_j\}$, again denoted by $\{f_j\}$, and for every $p\in\PC$
there is a sequence of functions $\{\vfi_{k,\,p}\}_{k=1}^\infty\subset%
\mbox{\rm Mon}(T;\RB^+)$ such that
  \begin{equation} \label{e:sta2}
\lim_{j\to\infty}V_{\vep_k,\,p}(f_j,T_t^-)=\vfi_{k,\,p}(t)\quad
\mbox{for all \,$t\in T$ and \,$k\in\NB$.}
  \end{equation}

We begin with $p=1$. The sequence
$\{t\mapsto V_{\vep_1,\,1}(f_j,T_t^-)\}_{j=j^*(\vep_1,1)}^\infty$ from
$\mbox{\rm Mon}(T;\RB^+)$  is uniformly bounded in $t\in T$ by $C(\vep_1,1)$,
and so, by the Helly selection principle, there are a subsequence
$\{J_{1,1}(j)\}_{j=1}^\infty$ of $\{j\}_{j=j^*(\vep_1,1)}^\infty$ and a function
$\vfi_{1,1}\in\mbox{\rm Mon}(T;\RB^+)$ such that
  $$\lim_{j\to\infty}V_{\vep_1,\,1}\bigl(f_{J_{1,1}(j)},T_t^-\bigr)=\vfi_{1,1}(t)\quad
     \mbox{for all \,$t\in T$.}$$
Now, choose the least number $j_{1,1}\in\NB$ such that $J_{1,1}(j_{1,1})%
\ge j^*(\vep_2,1)$.

\smallbreak
Inductively, assume that $k\ge2$ and a subsequence $\{J_{k-1,1}(j)\}_{j=1}^\infty$
of $\{j\}_{j=j^*(\vep_1,1)}^\infty$ and the number $j_{k-1,1}\in\NB$ such that
$J_{k-1,1}(j_{k-1,1})\ge j^*(\vep_k,1)$ have already been constructed. Since
$\{t\mapsto V_{\vep_k,\,1}(f_{J_{k-1,1}(j)},T_t^-)\}_{j=j_{k-1,1}}^\infty$ is a
sequence from $\mbox{\rm Mon}(T;\RB^+)$, uniformly bounded in $t\in T$ by
$C(\vep_k,1)$, the Helly selection principle implies the existence of a subsequence
$\{J_{k,1}(j)\}_{j=1}^\infty$ of $\{J_{k-1,1}(j)\}_{j=j_{k-1,1}}^\infty$ and a function
$\vfi_{k,1}\in\mbox{\rm Mon}(T;\RB^+)$ such that
  $$\lim_{j\to\infty}V_{\vep_k,\,1}\bigl(f_{J_{k,1}(j)},T_t^-\bigr)=\vfi_{k,1}(t)\quad
     \mbox{for all \,$t\in T$.}$$
Given $k\in\NB$, $\{J_{j,1}(j)\}_{j=k}^\infty$ is a subsequence of
$\{J_{k,1}(j)\}_{j=1}^\infty$, and so, the diagonal subsequence
$\{f_j^{(1)}\}_{j=1}^\infty\equiv\{f_{J_{j,1}(j)}\}_{j=1}^\infty$ of the original
sequence $\{f_j\}_{j=1}^\infty$ satisfies the condition
  \begin{equation} \label{e:ekek}
\lim_{j\to\infty}V_{\vep_k,\,1}\bigl(f_j^{(1)},T_t^-\bigr)=\vfi_{k,1}(t)\quad
     \mbox{for all \,$t\in T$ and \,$k\in\NB$.}
  \end{equation}

Taking into account \eq{e:ekek}, again inductively, assume that $p\in\PC$, $p\ge2$,
and a subsequence $\{f_j^{(p-1)}\}_{j=1}^\infty$ of $\{f_j^{(1)}\}_{j=1}^\infty$
(and, hence, of $\{f_j\}$) satisfying
  \begin{equation*}
\lim_{j\to\infty}V_{\vep_k,\,p-1}\bigl(f_j^{(p-1)},T_t^-\bigr)=\vfi_{k,p-1}(t)\quad
     \mbox{for all \,$t\in T$ and \,$k\in\NB$}
  \end{equation*}
is already chosen. The sequence
$\{t\mapsto V_{\vep_1,\,p}(f_j^{(p-1)},T_t^-)\}_{j=j^*(\vep_1,p)}^\infty$ of
nondecreasing functions on $T$ is uniformly bounded in $t\in T$ by $C(\vep_1,p)$,
and so, by the Helly selection principle, there are a subsequence
$\{f_{J_{1,p}(j)}^{(p-1)}\}_{j=1}^\infty$ of $\{f_j^{(p-1)}\}_{j=j^*(\vep_1,p)}^\infty$
(where $J_{1,p}:\NB\to\NB$ is a strictly increasing subsequence of
$\{j\}_{j=j^*(\vep_1,p)}^\infty$) and a function 
$\vfi_{1,p}\in\mbox{\rm Mon}(T;\RB^+)$ such that
  $$\lim_{j\to\infty}V_{\vep_1,\,p}\bigl(f_{J_{1,p}(j)}^{(p-1)},T_t^-\bigr)=\vfi_{1,p}(t)
     \quad\mbox{for all \,$t\in T$.}$$
Pick the least number $j_{1,p}\in\NB$ such that $J_{1,p}(j_{1,p})%
\ge j^*(\vep_2,p)$. Inductively, assume now that $k\ge2$ and a subsequence
$\{J_{k-1,p}(j)\}_{j=1}^\infty$ of $\{j\}_{j=j^*(\vep_1,p)}^\infty$ and the (least)
number $j_{k-1,p}\in\NB$ such that $J_{k-1,p}(j_{k-1,p})\ge j^*(\vep_k,p)$ are
already chosen. Since
$\{t\mapsto V_{\vep_k,\,p}(f_{J_{k-1,p}(j)}^{(p-1)},T_t^-)\}_{j=j_{k-1,p}}^\infty$ is a
sequence from $\mbox{\rm Mon}(T;\RB^+)$, uniformly bounded in $t\in T$ by
$C(\vep_k,p)$, the Helly selection principle implies the existence of a subsequence
$\{J_{k,p}(j)\}_{j=1}^\infty$ of $\{J_{k-1,p}(j)\}_{j=j_{k-1,p}}^\infty$ and a function
$\vfi_{k,p}\in\mbox{\rm Mon}(T;\RB^+)$ such that
  $$\lim_{j\to\infty}V_{\vep_k,\,p}\bigl(f_{J_{k,p}(j)}^{(p-1)},T_t^-\bigr)=\vfi_{k,p}(t)
     \quad\mbox{for all \,$t\in T$.}$$
Given $k\in\NB$, $\{J_{j,p}(j)\}_{j=k}^\infty$ is a subsequence of
$\{J_{k,p}(j)\}_{j=1}^\infty$, and so, the diagonal subsequence
$\{f_j^{(p)}\}_{j=1}^\infty\equiv\{f_{J_{j,p}(j)}^{(p-1)}\}_{j=1}^\infty$ of
$\{f_{J_{k,p}(j)}^{(p-1)}\}_{j=1}^\infty$ satisfies the equality
  \begin{equation*}
\lim_{j\to\infty}V_{\vep_k,\,p}\bigl(f_j^{(p)},T_t^-\bigr)=\vfi_{k,p}(t)\quad
     \mbox{for all \,$t\in T$ and \,$k\in\NB$.}
  \end{equation*}

Finally, the diagonal sequence $\{f_j^{(j)}\}_{j=1}^\infty$, again denoted
by $\{f_j\}\equiv\{f_j\}_{j=1}^\infty$, satisfies \eq{e:sta2}. 

\smallbreak
2. Let $Q$ be an at most countable dense subset of $T$. Note that any point $t\in T$,
which is not a limit point for $T$ (i.e., $T\cap(t-\delta,t+\delta)=\{t\}$ for some
$\delta>0$), belongs to~$Q$. Since, for any $k\in\NB$ and $p\in\PC$,
$\vfi_{k,p}\in\mbox{\rm Mon}(T;\RB^+)$, the set $Q_{k,p}\subset T$ of points of
discontinuity of $\vfi_{k,p}$ is at most countable. Setting
$S=Q\cup\bigcup_{k\in\NB}\bigcup_{p\in\PC}Q_{k,p}$, we find that $S$ is an  at most
countable dense subset of~$T$. Furthermore, if $S\ne T$, then every point
$t\in T\setminus S$ is a limit point for $T$ and
  \begin{equation} \label{e:dese}
\mbox{$\vfi_{k,p}$ is continuous on $T\setminus S$ for all $k\in\NB$ and $p\in\PC$.}
  \end{equation}
Since $S$ is at most countable and $\{f_j(s):j\in\NB\}$ is a relatively sequentially
compact subset of $X$ for all $s\in S$, with no loss of generality we may assume
(applying the Cantor diagonal process and passing to a subsequence of
$\{f_j\}$ if necessary) that, for each $s\in S$, $f_j(s)$ converges in $X$ as $j\to\infty$
to a unique point denoted by $f(s)\in X$, so that $f:S\to X$ and
$\lim_{j\to\infty}d_p(f_j(s),f(s))=0$ for all $p\in\PC$.

\smallbreak
If $S=T$, we turn to step~4 below which completes the proof.

\smallbreak
3. Assume that $S\ne T$ and $t\in T\setminus S$ is arbitrary. Let us prove that
$\{f_j(t)\}_{j=1}^\infty$ is a Cauchy sequence in $X$, i.e.,
  $$\lim_{j,j'\to\infty}d_p(f_j(t),f_{j'}(t))=0\quad\mbox{for \,all \,\,$p\in\PC$.}$$
Let us fix $p\in\PC$. Since $\vep_k\to0$ as $k\to\infty$, given $\eta>0$, choose and
fix $k=k(\eta)\in\NB$ such that $\vep_k\le\eta$. By \eq{e:dese}, $\vfi_{k,p}$ is
continuous at $t$, and so, by the density of $S$ in $T$, there is $s=s(\eta,k,t,p)\in S$
such that
  $$|\vfi_{k,p}(t)-\vfi_{k,p}(s)|\le\eta.$$
Property \eq{e:sta2} implies the existence of $j^1=j^1(\eta,k,t,s,p)\in\NB$ such that,
for all $j\ge j^1$, we have
  \begin{equation} \label{e:alf5}
|V_{\vep_k,\,p}(f_j,T_t^-)-\vfi_{k,p}(t)|\le\eta\quad\mbox{and}\quad
|V_{\vep_k,\,p}(f_j,T_s^-)-\vfi_{k,p}(s)|\le\eta.
  \end{equation}
Supposing (with no loss of generality) that $s<t$ and applying
Lemma~\ref{l:ess}\,(e) (in which $T$ is replaced by $T_t^-$, so that
$(T_t^-)_s^-=T_s^-$ and $(T_t^-)_s^+=T\cap[s,t]$), we get
  \begin{align*}
V_{\vep_k,\,p}(f_j,T\cap[s,t])&\le V_{\vep_k,\,p}(f_j,T_t^-)-
  V_{\vep_k,\,p}(f_j,T_s^-)\\[2pt]
&\le|V_{\vep_k,\,p}(f_j,T_t^-)-\vfi_{k,p}(t)|+|\vfi_{k,p}(t)-\vfi_{k,p}(s)|\\[2pt]
&\quad\,\,+|\vfi_{k,p}(s)-V_{\vep_k,\,p}(f_j,T_s^-)|\\[2pt]
&\le\eta+\eta+\eta=3\eta\quad\mbox{for all \,$j\ge j^1$.}
  \end{align*}
By the definition of $V_{\vep_k,\,p}(f_j,T\cap[s,t])$, for each $j\ge j^1$ there is a
function $g_j\in\mbox{\rm BV}_p(T\cap[s,t];X)$, also depending on $\eta$, $k$, $t$,
$s$ and $p$, such that
 $$d_{T\cap[s,t],\,p}(f_j,g_j)\le\vep_k\quad\mbox{and}\quad
    V_p(g_j,T\cap[s,t])\le V_{\vep_k,\,p}(f_j,T\cap[s,t])+\eta.$$
These inequalities, \eq{e:2-2} and the definition of $V_p(\cdot,\cdot)$ imply,
for all $j\ge j^1$,
  \begin{align}
d_p(f_j(s),f_j(t))&\le d_p(g_j(s),g_j(t))+2d_{T\cap[s,t],\,p}(f_j,g_j)\nonumber\\[2pt]
&\le V_p(g_j,T\cap[s,t])+2\vep_k\le(3\eta+\eta)+2\eta=6\eta. \label{e:6eta}
  \end{align}
Being convergent in the uniform space $X$, the sequence $\{f_j(s)\}_{j=1}^\infty$ is
Cauchy (cf.\ \cite[Theorem 6.21]{Kelley}), and so, there is $j^2=j^2(\eta,s,p)\in\NB$
such that $d_p(f_j(s),f_{j'}(s))\le\eta$ for all $j,j'\ge j^2$. The number
$j^3=\max\{j^1,j^2\}$ depends only on $\eta$, $t$ and $p$, and we find,
by the triangle inequality for $d_p$,
  \begin{align*}
d_p(f_j(t),f_{j'}(t))&\le d_p(f_j(t),f_j(s))+d_p(f_j(s),f_{j'}(s))+d_p(f_{j'}(s),f_{j'}(t))
  \\[2pt]
&\le 6\eta+\eta+6\eta=13\eta\quad\mbox{for all \,$j,j'\ge j^3$.}
  \end{align*}
Due to the arbitrariness of $p\in\PC$, this proves that $\{f_j(t)\}_{j=1}^\infty$ is a
Cauchy sequence in $X$. Taking into account that the set $\{f_j(t):j\in\NB\}$ is
relatively sequentially compact in $X$, we conclude that the sequence
$\{f_j(t)\}_{j=1}^\infty$ has a limit point in $X$, which we denote by $f(t)\in X$.
By \cite[Theorem 6.21]{Kelley}, a Cauchy sequence in a uniform space converges
to its limit point, and so, $\lim_{j\to\infty}d_p(f_j(t),f(t))=0$ for all $p\in\PC$.

\smallbreak
4. Since the uniform space $X$ is Hausdorff, we have shown at the end of steps 2
and~3 that the function $f:T=S\cup(T\setminus S)\to X$ is well-defined, and it is the
pointwise  limit on $T$ of a subsequence $\{f_{j_k}\}_{k=1}^\infty$ of the original
sequence $\{f_j\}_{j=1}^\infty$. It follows from Lemma~\ref{l:ess}(d) that,
given $p\in\PC$ and $\vep_0>0$,
  \begin{align*}
|f(T)|_p&\le\liminf_{k\to\infty}|f_{j_k}(T)|_p\le
  \liminf_{k\to\infty}V_{\vep_{{}_0},\,p}(f_{j_k},T)+2\vep_0\\[2pt]
&\le\limsup_{j\to\infty}V_{\vep_{{}_0},\,p}(f_j,T)+2\vep_0<\infty,
  \end{align*}
and so, $f$ is a bounded function on $T$ (i.e., $f\in\mbox{\rm B}(T;X)$).

\smallbreak
Now, we prove that $f$ is {\em regulated\/} on $T$. Given $\tau\in T$, which is a left
limit point for $T$, let us show that $d_p(f(s),f(t))\to0$ as $T\ni s,t\to\tau-0$ for all
$p\in\PC$ (similar arguments apply in the case when $\tau'\in T$ is a right limit point
for~$T$). Given $p\in\PC$, this is equivalent to showing that for every $\eta>0$
there is $\delta=\delta(\eta,p)>0$ such that $d_p(f(s),f(t))\le7\eta$ for all
$s,t\in T\cap(\tau-\delta,\tau)$ with $s<t$. Recall that the (finally) extracted
subsequence of the original sequence $\{f_j\}$, again denoted by $\{f_j\}$ here,
satisfies condition \eq{e:sta2}, and $f_j\to f$ pointwise on~$T$.

\smallbreak
Let $p\in\PC$ and $\eta>0$ be arbitrarily fixed. Since $\vep_k\to0$ as $k\to\infty$,
pick and fix natural $k=k(\eta)$ such that $\vep_k\le\eta$, and since
$\vfi_{k,p}\in\mbox{\rm Mon}(T;\RB^+)$ and $\tau\in T$ is a left limit point for $T$,
the left limit $\lim_{T\ni t\to\tau-0}\vfi_{k,p}(t)\in\RB^+$ exists. Hence, there is
$\delta=\delta(\eta,k,p)>0$ such that $|\vfi_{k,p}(t)-\vfi_{k,p}(s)|\le\eta$ for all
$s,t\in T\cap(\tau-\delta,\tau)$. By virtue of \eq{e:sta2}, for any
$s,t\in T\cap(\tau-\delta,\tau)$ there is $j^1=j^1(\eta,k,s,t,p)\in\NB$ such that,
if $j\ge j^1$, then the inequalities \eq{e:alf5} hold. Arguing exactly the same way as
in between the lines \eq{e:alf5} and \eq{e:6eta}, we find that
$d_p(f_j(s),f_j(t))\le6\eta$ for all $j\ge j^1$. Noting that $f_j(s)\to f(s)$ and
$f_j(t)\to f(t)$ in $X$ as $j\to\infty$, by the triangle inequality for $d_p$, we get
  $$|d_p(f_j(s),f_j(t))-d_p(f(s),f(t))|\le d_p(f_j(s),f(s))+d_p(f_j(t),f(t))\to0$$
as $j\to\infty$. So, there is $j^2=j^2(\eta,s,t,p)\in\NB$ such that
$d_p(f(s),f(t))\le d_p(f_j(s),f_j(t))+\eta$ for all $j\ge j^2$. Thus, choosing
$j\ge\max\{j^1,j^2\}$, we obtain $d_p(f(s),f(t))\le6\eta+\eta=7\eta$.

\smallbreak
This completes the proof of Theorem~\ref{th:2}.
\qed\end{pthscnd}

\begin{rem}
Theorem \ref{th:2} is an extension of the Helly-type selection principles from
\cite[Theorem 3.8]{Fr} (for $T=[a,b]$ and $X=\RB^N$) and \cite[Theorem~3]{Studia}
(for any $T\subset\RB$ and a metric space $X$).
It also extends Theorem~8 from \cite{MatTr}, in which,
under the assumptions of Theorem~\ref{th:2}, condition \eq{e:Vepfj} is replaced
by a more stringent one: $C_p\equiv\sup_{j\in\NB}V_p(f_j,T)<\infty$ for all $p\in\PC$.
In fact, by Lemma~\ref{l:ess}\,(c), $V_{\vep,\,p}(f_j,T)\le C_p$ for all $j\in\NB$,
$\vep>0$ and $p\in\PC$, and so, \eq{e:Vepfj} is fulfilled. Now, if, according to
Theorem~\ref{th:2}, a subsequence of $\{f_j\}$ converges pointwise on $T$
to a function $f\in X^T$, then property \eq{e:Vplsc} implies $V_p(f,T)\le C_p$
for all $p\in\PC$, i.e., $f\in\mbox{\rm BV}(T;X)$.
\end{rem}

\begin{rem}
Condition \eq{e:Vepfj} in Theorem \ref{th:2} is {\em necessary\/} for the
{\em uniformly\/} convergent sequence $\{f_j\}\subset X^T$ in the following sense:
if $f_j\rra f$ on $T$, where $f\in X^T$ is such that $V_{\vep,\,p}(f,T)<\infty$ for all
$\vep>0$ and $p\in\PC$, then, by virtue of Lemma~\ref{l:unif}\,(a), for all
$0<\vep'<\vep$ and $p\in\PC$, we have
  $$\limsup_{j\to\infty}V_{\vep,\,p}(f_j,T)\le V_{\vep-0,\,p}(f,T)
     \le V_{\vep'\!,\,p}(f,T)<\infty.$$
\end{rem}

\begin{exa}
Condition \eq{e:Vepfj} in Theorem \ref{th:2} is {\em not necessary\/} for the
{\em pointwise\/} convergent sequence $\{f_j\}\subset X^T$ (cf.\ also Examples
4.7, 4.10 and 4.11 in \cite{arXiv19}). This can be illustrated by the sequence $\{f_j\}$
from Example~\ref{exa4}, where $I=[0,1]$. We assert that if $p\in\PC$ is such that
$d_p(x,y)>0$, and $0<\vep<d_p(x,y)/2$, then
$\lim_{j\to\infty}V_{\vep,\,p}(f_j,I)=\infty$. In fact, given $j\in\NB$, we set
$t_k=k/j!$ (so that $f_j(t_k)=x$) for $k=0,1,\dots,j!$, and
$s_k=(t_{k-1}+t_k)/2=(k-\frac12)/j!$ (so that $f_j(s_k)=y$) for $k=1,2,\dots,j!$.
This induces a partition of $I=[0,1]$:
  $$0=t_0<s_1<t_1<s_2<t_2<\ldots<s_{j!-1}<t_{j!-1}<s_{j!}<t_{j!}=1.$$
Supposing $g\in X^I$ is arbitrary such that $d_{I,\,p}(f_j,g)\le\vep$, by virtue
of the definition of $V_p(g,I)$ and \eq{e:2-2}, we find
  $$V_p(g,I)\ge\sum_{k=1}^{j!}d_p(g(t_k),g(s_k))\ge\sum_{k=1}^{j!}
     \bigl(d_p(f(t_k),f(s_k))-2\vep\bigr)\!=\!j!\cdot(d_p(x,y)-2\vep).$$
By definition \eq{d:Vep}, $V_{\vep,\,p}(f_j,I)\ge j!\cdot(d_p(x,y)-2\vep)$, which
proves our assertion.
\end{exa}

\begin{exa}
Theorem \ref{th:2} is inapplicable to the sequence $\{f_j\}$ from Example~\ref{ex:vez},
because (although $f_j\rra f\equiv\DC_{x,y}$ on $I$)
$\lim_{j\to\infty}V_\vep(f_j,I)=\infty$ for $\vep=\|x-y\|/2$.
 The reason here is that the limit function
$\DC_{x,y}$ is not regulated if $x\ne y$.

\smallbreak
Nevertheless, Theorem \ref{th:2} can be applied to sequences of {\em non-regulated\/}
functions. Let $\{x_j\}$ and $\{y_j\}$ be two sequences in the uniform space $X$
such that $x_j\ne y_j$ for all $j\in\NB$, $x_j\to x$ and $y_j\to x$ in $X$ as
$j\to\infty$ for some $x\in X$. The sequence $f_j=\DC_{x_j,y_j}$ converges
uniformly on $I=[a,b]$ to the constant function $f(t)\equiv x$ on $I$:
  $$d_{I,\,p}(f_j,f)=\max\{d_p(x_j,x),d_p(y_j,x)\}\to0\quad\!
    \mbox{as $j\to\infty$ for all \,$p\in\PC$.}$$
Given $p\in\PC$ and $\vep>0$, there is $j_0=j_0(p,\vep)\in\NB$ such that
$d_p(x_j,y_j)\le\vep$ for all $j\ge j_0$, and so, by \eq{e:inf0}, $V_{\vep,\,p}(f_j,I)=0$
for all $j\ge j_0$. This implies condition \eq{e:Vepfj}:
  $$\limsup_{j\to\infty}V_{\vep,\,p}(f_j,I)\le\sup_{j\ge j_0}V_{\vep,\,p}(f_j,I)=0.$$
\end{exa}

Applying Theorem \ref{th:2} and the diagonal process over expanding intervals,
we get the following {\em local\/} version of Theorem~\ref{th:2}:

\begin{thm}
Under the assumptions of Theorem~{\rm\ref{th:2}}, suppose that condition\/
\eq{e:Vepfj} is replaced by\/ {\rm(}the local one{\rm)}
  $$\limsup_{j\to\infty}V_{\vep,\,p}(f_j,T\cap[a,b])<\infty\,\,\,\mbox{for all
     $a,b\in T$, $a\le b$, $\vep>0$ and $p\in\PC$.}$$
Then, there is a subsequence of $\{f_j\}$ which converges pointwise on $T$ to a function
$f\in\mbox{\rm Reg}(T;X)$ such that $f\in\mbox{\rm B}(T\cap[a,b];X)$
for all $a,b\in T$, $a\le b$.
\end{thm}



\bigbreak

\end{document}